\documentclass[12pt,reqno]{amsart}
\usepackage{enumerate, latexsym, amsmath, amsfonts, amssymb, amsthm, graphicx, color}
\def\pmod #1{\ ({\rm{mod}}\ #1)}
\def\Z{\Bbb Z}
\def\N{\Bbb N}

\def\l{\left}
\def\r{\right}
\def\bg{\bigg}
\def\({\bg(}
\def\){\bg)}

\def\f{\frac}

\def\ls{\leqslant}
\def\gs{\geqslant}
\def\se {\subseteq}
\def\sm{\setminus}
\def\bi{\binom}
\def\al{\alpha}

\def\eq{\equiv}

\def\da{\delta}

\def\Proof{\noindent{\it Proof}}

\theoremstyle{plain}
\newtheorem{theorem}{Theorem}

\newtheorem{lemma}{Lemma}
\newtheorem{corollary}{Corollary}

\theoremstyle{definition}

\theoremstyle{remark}
\newtheorem{remark}{Remark}

\makeatletter
\@namedef{subjclassname@2010}{%
  \textup{2010} Mathematics Subject Classification}
\makeatother
 \vspace{4mm}

\begin{document}
 \baselineskip=17pt
\hbox{Acta Arith. 180(2017), no.\,3, 229--249.}
\medskip

\title[Sums of four polygonal numbers with coefficients]{Sums of four polygonal numbers with coefficients}

\author[X.-Z. Meng]{Xiang-Zi Meng}
\address{(Xiang-Zi Meng) Department of Mathematics\\ Nanjing University\\
Nanjing 210093, People's Republic of China}
\email{xzmeng@smail.nju.edu.cn}

\author[Z.-W. Sun]{Zhi-Wei Sun}
\address{(Zhi-Wei Sun, corresponding author) Department of Mathematics\\ Nanjing University\\
Nanjing 210093, People's Republic of China}
\email{zwsun@nju.edu.cn}

\date{}

\begin{abstract}
Let $m\gs3$ be an integer. The polygonal numbers of order $m+2$ are given by $p_{m+2}(n)=m\binom n2+n$ $(n=0,1,2,\ldots)$.
 A famous claim of Fermat proved by Cauchy asserts that each nonnegative integer is the sum of $m+2$ polygonal numbers of order $m+2$.
 For $(a,b)=(1,1),(2,2),(1,3),(2,4)$, we study whether any sufficiently large integer can be expressed as
$$p_{m+2}(x_1)+p_{m+2}(x_2)+ap_{m+2}(x_3)+bp_{m+2}(x_4)$$
with $x_1,x_2,x_3,x_4$ nonnegative integers.
We show that the answer is positive if $(a,b)\in\{(1,3),(2,4)\}$, or $(a,b)=(1,1)\ \&\ 4\mid m$, or $(a,b)=(2,2)\ \&\ m\not\eq2\pmod4$.
In particular, we confirm a conjecture of Z.-W. Sun which states that any natural number can be written as
$p_6(x_1)+p_6(x_2)+2p_6(x_3)+4p_6(x_4)$ with $x_1,x_2,x_3,x_4$ nonnegative integers.
\end{abstract}



\maketitle

\section{Introduction}
\setcounter{lemma}{0}
\setcounter{theorem}{0}
\setcounter{corollary}{0}
\setcounter{remark}{0}
\setcounter{equation}{0}
\setcounter{conjecture}{0}

Let $m\in\Z^+=\{1,2,3,\ldots\}$.
The {\it polygonal numbers of order $m+2$} (or {\it $(m+2)$-gonal numbers}), which are constructed geometrically from the regular polygons with $m+2$ sides,
are given by
\begin{equation}\label{1.1}p_{m+2}(n):=m\binom n2+n=\frac{mn^2-(m-2)n}2\ \ \mbox{for}\ n\in\N=\{0,1,2,\ldots\}.\end{equation}
Clearly,
\[p_{m+2}(0)=0,\ p_{m+2}(1)=1,\ p_{m+2}(2)=m+2,\ p_{m+2}(3)=3m+3,\]
and $p_{m+2}(x)$ with $x\in\Z$ are called {\it generalized $(m+2)$-gonal numbers}.
It is easy to see that generalized hexagonal numbers coincide with triangular numbers (i.e., those $p_3(n)=n(n+1)/2$ with $n\in\N$). Note that
\[p_4(n)=n^2,\ p_5(n)=\frac {n(3n-1)}2,\ p_6(n)=n(2n-1)=p_3(2n-1).\]
Fermat's claim that each $n\in\N$ can be written as the sum of $m+2$ polygonal numbers of order $m+2$ was proved by Lagrange in the case $m=2$, by Gauss in the case $m=1$, and by Cauchy in the case $m\gs3$ (cf. \cite[pp.\,3-35]{N96} and \cite[pp. 54-57]{MW}). In 1830 Legendre refined Cauchy's polygonal number theorem by showing that any integer $N\gs 28m^3$ with $m\gs 3$
can be written as $$p_{m+2}(x_1)+p_{m+2}(x_2)+p_{m+2}(x_3)+p_{m+2}(x_4)+\da_m(N)$$
where $x_1,x_2,x_3,x_4\in\N$, $\da_m(N)=0$ if $2\nmid m$, and $\da_m(N)\in\{0,1\}$ if $2\mid m$.
Nathanson (\cite{N87} and \cite[p.\,33]{N96}) simplified the proofs of Cauchy's and Legendre's theorems.

In 1917 Ramanujan \cite{R} listed 55 possible quadruples $(a,b,c,d)$ of positive integers with $a\ls b\ls c\ls d$
such that any $n\in\N$ can be written as $ax^2+by^2+cz^2+dw^2$ with $x,y,z,w\in\Z$, and
54 of them were later confirmed by Dickson \cite{D27} while the remaining one on the list was actually wrong.

Recently, Sun \cite{S16} showed that any positive integer can be written as the sum of four generalized octagonal numbers one of which is odd. He also proved that for many triples
$(b,c,d)$ of positive integers (including (1,1,3), (1,2,2) and (1,2,4)) we have
$$\{p_8(x_1)+bp_8(x_2)+cp_8(x_3)+dp_8(x_4):\ x_1,x_2,x_3,x_4\in\Z\}=\N.$$
In \cite[Conjecture 5.3]{S16}, Sun conjectured that any $n\in\N$ can be written as $p_6(x_1)+p_6(x_2)+2p_6(x_3)+4p_6(x_4)$ with $x_1,x_2,x_3,x_4\in\N$.

Motivated by the above work, for $(a,b)=(1,1),(2,2),(1,3),(2,4)$ and $m\in\{3,4,5,\ldots\}$, we study whether any sufficiently large integer can be written as
$$p_{m+2}(x_1)+p_{m+2}(x_2)+ap_{m+2}(x_3)+bp_{m+2}(x_4)\ \ \mbox{with}\ x_1,x_2,x_3,x_4\in \N.$$

Now we state our main results.

\begin{theorem} \label{Th1.1} Let $m\in\Z^+$ with $4\mid m$.

{\rm (i)} Any integer $N\gs28m^3$ can be expressed as
\begin{equation}\label{1.2}p_{m+2}(x_1)+p_{m+2}(x_2)+p_{m+2}(x_3)+p_{m+2}(x_4)\ \ (x_1,x_2,x_3,x_4\in\N).\end{equation}

{\rm (ii)} There are infinitely many positive integers not of the form
$p_{m+4}(x_1)+p_{m+4}(x_2)+p_{m+4}(x_3)+p_{m+4}(x_4)$ with $x_1,x_2,x_3,x_4\in\N$.
\end{theorem}
\begin{remark}\label{Re1.1} This can be viewed as a supplement to Legendre's theorem. By Theorem 1.1(ii), there are infinitely many positive integers
none of which is the sum of four octagonal numbers; in contrast, Sun \cite{S16} showed that any $n\in\N$ is the sum of four generalized octagonal numbers.
\end{remark}

\begin{corollary}\label{Cor1.1} We have
\begin{equation}\label{1.3}\begin{aligned}&\{p_6(x_1)+p_6(x_2)+p_6(x_3)+p_6(x_4):\ x_1,x_2,x_3,x_4\in\N\}
\\=&\N\sm\{5,\,10,\,11,\,20,\,25,\,26,\,38,\,39,\,54,\,65,\,70,\,114,\,130\}
\end{aligned}\end{equation}
and hence any $n\in\N$ can be written as the sum of a triangular number and three hexagonal numbers.
Also, any integer $n>2146$ can be written as the sum of four decagonal numbers and thus
\begin{equation}\label{1.4}
\{p_{10}(x_1)+p_{10}(x_2)+p_{10}(x_3)+p_{10}(x_4):\ x_1,x_2,x_3,x_4\in\Z\}=\N\sm\{5,\,6,\,26\}.
\end{equation}
\end{corollary}
\Proof. Via a computer, we can easily verify that
$$5,\,10,\,11,\,20,\,25,\,26,\,38,\,39,\,54,\,65,\,70,\,114,\,130$$
are the only natural numbers smaller than $28\times4^3$ which
cannot be written as the sum of four hexagonal numbers, but all these numbers can be expressed as the sum of a triangular number and three hexagonal numbers. Also, every $n=2147,\ldots,28\times 8^3-1$ is the sum of four decagonal numbers,
and 5, 6 and 26 are the only natural numbers smaller than 2147 which cannot be written as the sum of four generalized decagonal numbers.
Now it suffices to apply Theorem 1.1 with $m=4,8$. \qed

\begin{remark}\label{Rem1.2} Sun \cite[Conjecture 1.10]{S15} conjectured that any $n\in\N$ can be written as the sum of two triangular numbers and a hexagonal number.
Krachun \cite{K} proved that
\begin{align*}&\{p_6(-w)+p_6(-x)+p_6(y)+p_6(z):\ w,x,y,z\in\N\}
\\=&\{p_6(-w)+2p_6(-x)+p_6(y)+2p_6(z):\ w,x,y,z\in\N\}=\N,
\end{align*}
which was first conjectured by the second author \cite{S16}.
\end{remark}

\begin{theorem} \label{Th1.2} Let $m\gs 3$ be an integer.

{\rm (i)} Suppose that $2\nmid m$ or $4\mid m$. Then any integer $N\gs1628m^3$ can be written as
\begin{equation}\label{1.5}p_{m+2}(x_1)+p_{m+2}(x_2)+2p_{m+2}(x_3)+2p_{m+2}(x_4)\ \ (x_1,x_2,x_3,x_4\in\N).\end{equation}

{\rm (ii)} If $m\eq2\pmod 4$, then there are infinitely many positive integers not represented by
$p_{m+2}(x_1)+p_{m+2}(x_2)+2p_{m+2}(x_3)+2p_{m+2}(x_4)$ with $x_1,x_2,x_3,x_4\in\N$.
\end{theorem}
\begin{remark}\label{Rem1.3} Actually our proof of Theorem 1.2(i) given in Section 3
allows us to replace $1628m^3$ by $418m^3$ in the case $m\eq1\pmod2$. By Theorem 1.2(ii), there are infinitely many positive
integers not represented by $p_8(x_1)+p_8(x_2)+2p_8(x_3)+2p_8(x_4)$ with $x_1,x_2,x_3,x_4\in\N$; in contrast, Sun \cite{S16}
proved that any $n\in\N$ can be written as $p_8(x_1)+p_8(x_2)+2p_8(x_3)+2p_8(x_4)$ with $x_1,x_2,x_3,x_4\in\Z$.
\end{remark}

\begin{corollary}\label{Cor1.2} We have
\begin{equation}\label{1.6}\{p_5(x_1)+p_5(x_2)+2p_5(x_3)+2p_5(x_4):\ x_1,x_2,x_3,x_4\in\N\}=\N,
\end{equation}
\begin{equation}\label{1.7}
\begin{aligned}&\{p_6(x_1)+p_6(x_2)+2p_6(x_3)+2p_6(x_4):\ x_1,x_2,x_3,x_4\in\N\}
\\&\qquad=\N\sm\{22,\,82,\,100\},
\end{aligned}\end{equation}
and
\begin{equation}\label{1.8}
\begin{aligned}&\{p_7(x_1)+p_7(x_2)+2p_7(x_3)+2p_7(x_4):\ x_1,x_2,x_3,x_4\in\N\}
\\&\ \ \ \ =\N\sm\{13,\,26,\,31,\,65,\,67,\,173,\,175,\,215,\,247\}.
\end{aligned}
\end{equation}
Also, for each $k=9,10,11$, any integer $n>C_k$ can be written as $p_k(x_1)+p_k(x_2)+2p_k(x_3)+2p_k(x_4)$ with $x_1,x_2,x_3,x_4\in\N$,
where $C_9=925$, $C_{10}=840$ and $C_{11}=1799$.
Therefore,
\begin{align}\label{1.9}\{p_3(w)+p_6(x)+2p_6(y)+2p_6(z):\ w,x,y,z\in\N\}=&\N,
\\\label{1.10}\{2p_3(w)+p_6(x)+p_6(y)+2p_6(z):\ w,x,y,z\in\N\}=&\N,
\\\label{1.11}\{p_7(w)+p_7(x)+2p_7(y)+2p_7(z):\ w\in\Z\ \&\ x,y,z\in\N\}=&\N,
\\\label{1.12}\{p_9(w)+2p_9(x)+p_9(y)+2p_9(z):\ w,x\in\Z\ \&\ y,z\in\N\}=&\N,
\\\label{1.13}\{p_{10}(w)+2p_{10}(x)+p_{10}(y)+2p_{10}(z):\ w,x\in\Z\ \&\ y,z\in\N\}=&\N,
\end{align}
and
\begin{equation}\label{1.14}\{p_{11}(w)+2p_{11}(x)+p_{11}(y)+2p_{11}(z):\ w,x\in\Z\ \&\ y,z\in\N\}=\N\sm\{7\}.\end{equation}
\end{corollary}
\Proof. Note that $\{p_6(w):\, w\in\Z\}=\{p_3(w):\, w\in\N\}$.
It suffices to apply Theorem 1.2(i) with $m\in\{3,4,5,7,8,9\}$ and check those $n\in\N$ with $n<1628m^3$ via a computer. \qed

\begin{remark}\label{Rem1.4} (1.6) appeared as part of \cite[Conjecture 5.2(ii)]{S16}, and it indicates that the set
$\{p_5(x)+2p_5(y):\ x,y\in\N\}$ is an additive base of order 2. For positive integers $a,b,c$ with
$\{ap_5(x)+bp_5(y)+cp_5(z):\ x,y,z\in\Z\}=\N$, see \cite{S15} and \cite{GS}.
\end{remark}

\begin{theorem} \label{Th1.3} Let $m\gs3$ be an integer. Then each integer $N\gs924m^3$ can be expressed as
\begin{equation}\label{1.15}p_{m+2}(x_1)+p_{m+2}(x_2)+p_{m+2}(x_3)+3p_{m+2}(x_4)\ \ (x_1,x_2,x_3,x_4\in\N).\end{equation}
\end{theorem}

\begin{corollary}\label{Cor1.3} We have
\begin{equation}\label{1.16}\{p_5(x_1)+p_5(x_2)+p_5(x_3)+3p_5(x_4):\ x_1,x_2,x_3,x_4\in\N\}=\N\sm\{19\},
\end{equation}
\begin{equation}\label{1.17}\begin{aligned}&\{p_6(x_1)+p_6(x_2)+p_6(x_3)+3p_6(x_4):\ x_1,x_2,x_3,x_4\in\N\}
\\&\ \ =\N\sm\{14,\,23,\,41,\,42,\,83\}\end{aligned}
\end{equation}
and
\begin{equation}\label{1.18}\begin{aligned}&\{p_7(x_1)+p_7(x_2)+p_7(x_3)+3p_7(x_4):\ x_1,x_2,x_3,x_4\in\N\}
\\=&\N\sm\{13,\,16,\, 27,\, 31,\,33,\,49,\,50,\,67,\,87,\,178,\,181,\,259\}.\end{aligned}
\end{equation}
Also, for each $k=8,9,10$, any integer $n>M_k$ can be written as $p_k(x_1)+p_k(x_2)+p_k(x_3)+3p_k(x_4)$ with $x_1,x_2,x_3,x_4\in\N$,
where $M_8=435,$ $M_9=695$ and $M_{10}=916$. Therefore
\begin{align}\label{1.19}\{p_7(w)+p_7(x)+p_7(y)+3p_7(z):\ w\in\Z\ \&\ x,y,z\in\N\}=&\N,
\\\label{1.20}\{p_9(x_1)+p_9(x_2)+p_9(x_3)+3p_9(x_4):\ x_1,x_2,x_3,x_4\in\Z\}=&\N\sm\{17\},
\\\label{1.21}\{p_{10}(x_1)+p_{10}(x_2)+p_{10}(x_3)+3p_{10}(x_4):\ x_1,x_2,x_3,x_4\in\Z\}=&\N\sm\{16,\,19\}.
\end{align}
\end{corollary}
\Proof. It suffices to apply Theorem 1.3 with $m\in\{3,4,5,6,7,8\}$ and check those $n\in\N$ with $n<924m^3$ via a computer. \qed

\begin{remark}\label{Rem1.5} Guy \cite{G} thought that $10$, $16$ and $76$ might be the only natural numbers which cannot be written
as the sum of three generalized heptagonal numbers. The second author \cite[Remark 5.2 and Conjecture 1.2]{S16} conjectured that
$\{p_7(x)+p_7(y)+p_7(z):\ x,y,z\in\Z\}=\N\sm\{10,16,76,307\}$
and $$\{p_8(x)+p_8(y)+3p_8(z):\ x,y,z\in\Z\}=\N\sm\{7,14,18,91\}.$$
\end{remark}

\begin{theorem} \label{Th1.4} Let $m\gs 3$ be an integer. Then any integer $N\gs1056m^3$ can be written as
\begin{equation}\label{1.22}p_{m+2}(x_1)+p_{m+2}(x_2)+2p_{m+2}(x_3)+4p_{m+2}(x_4)\ \ (x_1,x_2,x_3,x_4\in\N).\end{equation}
\end{theorem}

\begin{corollary}\label{Cor1.4} We have
\begin{align}\label{1.23}\{p_5(x_1)+p_5(x_2)+2p_5(x_3)+4p_5(x_4):\ x_1,x_2,x_3,x_4\in\N\}=&\N,
\\\label{1.24}\{p_6(x_1)+p_6(x_2)+2p_6(x_3)+4p_6(x_4):\ x_1,x_2,x_3,x_4\in\N\}=&\N,
\\\label{1.25}\{p_7(x_1)+p_7(x_2)+2p_7(x_3)+4p_7(x_4):\ x_1,x_2,x_3,x_4\in\N\}=&\N\sm\{17,51\},
\end{align}
\begin{equation}\label{1.26}\begin{aligned}&\{p_8(x_1)+p_8(x_2)+2p_8(x_3)+4p_8(x_4):\ x_1,x_2,x_3,x_4\in\N\}
\\&\quad \ \ \ =\N\sm\{19,\,30,\,39,\,59,\,78,\,91\},
\end{aligned}\end{equation}
and
\begin{equation}\label{1.27}\begin{aligned}&\{p_9(x_1)+p_9(x_2)+2p_9(x_3)+4p_9(x_4):\ x_1,x_2,x_3,x_4\in\N\}
\\=&\N\sm\{17,\,21,\,34,\,41,\,44,\,67,\,89,\,104,\,119,\,170,\,237,\,245,\,290\}.
\end{aligned}\end{equation}
Also, for each $k=10,11,12$ any integer $n>N_k$ can be written as $p_{k}(x_1)+p_{k}(x_2)+2p_{k}(x_3)+4p_{k}(x_4)$ with $x_1,x_2,x_3,x_4\in\N$,
where $N_{10}=333$, $N_{11}=734$ and $N_{12}=1334$. Therefore,
\begin{equation}\label{1.28}\{p_k(w)+p_k(x)+2p_k(y)+4p_k(z):\, w\in\Z\ \&\ x,y,z\in\N\}=\N\end{equation}
for $k=7,9$, and
\begin{equation}\label{1.29}\{p_k(w)+p_k(x)+2p_k(y)+4p_k(z):\ w,x,y,z\in\Z\}=\N\end{equation}
for $k=8,10,11,12$.
\end{corollary}
\Proof. It suffices to apply Theorem 1.4 with $m\in\{3,\ldots,10\}$ and check those $n\in\N$ with $n<1056m^3$ via a computer. \qed

\begin{remark}\label{Rem1.6} (\ref{1.23}) and (\ref{1.24}) were first conjectured by the second author \cite[Conjecture 5.2(ii) and Conjecture 5.3]{S16}.
Sun \cite[Remark 5.2]{S16} also conjectured that
$$\{p_7(x)+2p_7(y)+4p_7(z):\ x,y,z\in\Z\}=\N\sm\{131,146\}.$$
\end{remark}

We will show Theorems 1.1-1.4 in Sections 2-5 respectively.

Throughout this paper, for a prime $p$ and $a,n\in\N$, by $p^a\|n$ we mean $p^a\mid n$ and $p^{a+1}\nmid n$.
For example,  $4\|n$ if and only if $n\eq4\pmod8$.
\medskip

\maketitle

\section{Proof of Theorem 1.1}
\setcounter{lemma}{0}
\setcounter{theorem}{0}
\setcounter{corollary}{0}
\setcounter{remark}{0}
\setcounter{equation}{0}
\setcounter{conjecture}{0}

We first give a lemma which is a slight variant of \cite[Lemma 1.10]{N96}.

\begin{lemma}\label{Lem2.1} Let $l,m,N\in\Z^+$ with $N\gs7l^2m^3$. Then the length of the interval
\begin{equation}\label{2.1}I_1= \left[ \frac12+\sqrt{\frac{6N}m-3}, \ \frac23+\sqrt\frac {8N}m \right]
\end{equation}
is greater than $lm$.
\end{lemma}
\Proof. Let $L_1$ denote the length of the interval $I_1$. Then $L_1=\sqrt{8x}-\sqrt{6x-3}+1/6$, where
$x=N/m\gs 7l^2m^2$. Let $l_0=lm-1/6$. Then
\begin{align*}L_1>lm&\iff \sqrt{8x}>\sqrt{6x-3}+l_0
\\&\iff 2x+3-l_0^2>2l_0\sqrt{6x-3}
\\&\iff 4x(x+3-7l_0^2)+(l_0^2-3)^2+12l_0^2>0
\end{align*}
As $x\gs 7l_0^2$, by the above we have $L_1>lm$.
This completes the proof. \qed

The following lemma is a slight modification of \cite[Lemma 1.11]{N96}.

\begin{lemma}\label{Lem2.2} Let $a,b,m,N\in\Z^+$ with $m\gs 3$ and
\begin{equation*}N=\frac m2(a-b)+b\gs\f23m.\end{equation*}
Suppose that $b$ belongs to the interval $I_1$ given by $(\ref{2.1})$.
Then
\begin{equation} \label{ineq}b^2<4a \quad\mbox{and}\quad 3a<b^2+2b+4. \end{equation}
\end{lemma}

\Proof. Observe that
\[a=\left(1-\frac2m\right)b+\frac {2N}m\]
and
\begin{align*}b^2-4a=&b^2-4\left(1-\frac2m\right)b-\frac{8N}m
\\=&\l(b-2\l(1-\f2m\r)\r)^2-4\l(\l(1-\f2m\r)^2+\f{2N}m\r).
\end{align*}
As $m\gs3$ and $b\in I_1$, we have
$$b\ls\frac23+\sqrt\frac {8N}m<2\left(1-\frac2m\right)+2\sqrt{\left(1-\frac2m\right)^2+\frac {2N}m}$$
and hence $b^2-4a<0$.
On the other hand, since $(1/2-3/m)^2<1$ and $b\in I_1$ we have
$$ b\gs\frac12+\sqrt{\frac{6N}m-3}>\frac12-\frac3m+\sqrt{\left(\frac12-\frac3m\right)^2-4+\frac{6N}m}$$
and hence
\begin{align*}b^2+2b+4-3a=&b^2-\left(1-\frac{6}m\right)b+\left(4-\frac{6N}m\right)
\\=&\l(b-\l(\f12-\f3m\r)\r)^2-\l(\l(\f12-\f3m\r)^2-4+\f{6N}m\r)
\\>&0.
\end{align*}
This proves (\ref{ineq}). \qed

\begin{lemma}\label{Lem2.3} Let $a,b,c$ be positive integers and let $x,y,z$ be real numbers. Then we have the inequality
\begin{equation}\label{2.3}(ax+by+cz)^2\ls(a+b+c)(ax^2+by^2+cz^2).
\end{equation}
\end{lemma}
\Proof. By the Cauchy-Schwarz inequality (cf. \cite[p.\,178]{N96}),
\begin{align*}&\l(\sqrt a(\sqrt ax)+\sqrt b(\sqrt by)+\sqrt c(\sqrt cz)\r)^2
\\\ls&\l((\sqrt a)^2+(\sqrt b)^2+(\sqrt c)^2\r)\l((\sqrt ax)^2+(\sqrt by)^2+(\sqrt cz)^2\r).
\end{align*}
This yields the desired (\ref{2.3}). \qed

The following lemma with $2\nmid ab$ is usually called Cauchy's Lemma (cf. \cite[pp. 31--34]{N96}).

\begin{lemma}\label{Lem2.4} Let $a$ and $b$ be positive integers satisfying $(\ref{ineq})$.
Suppose that $2\nmid ab$, or $2\|a$ and $2\mid b$.
Then there exist $s,t,u,v\in\N$ such that
\begin{equation}\label{2.4}a=s^2+t^2+u^2+v^2\ \ \mbox{and}\ \ b=s+t+u+v.\end{equation}
\end{lemma}
\Proof. By the Gauss-Legendre theorem (cf. \cite[Section 1.5]{N96}), we have
\begin{equation}\label{2.5}\{x^2+y^2+z^2:\ x,y,z\in\Z\}=\N\setminus\{4^k(8l+7):\ k,l\in\N\}.\end{equation}
We claim that there are $x,y,z\in\Z$ with $4a-b^2=x^2+y^2+z^2$ such that all the numbers
\begin{equation}\label{2.6}
\begin{aligned}&s=\frac{b+ x+ y+ z}4,\ \ t=\frac{b+ x- y-z}4,
\\& u=\frac{b- x+y- z}4,\ \  v=\frac{b- x-y+ z}4
\end{aligned}
\end{equation}
are integers.
\medskip

{\it Case} 1. $2\nmid ab$.

In this case, $4a-b^2\equiv3\pmod{8}$ and hence $4a-b^2=x^2+y^2+z^2$ for some $x,y,z\in\Z$ with $2\nmid xyz$.
Without loss of generality, we may assume that $x\eq y\eq z\eq b\pmod 4$. (If $x\eq-b\pmod 4$ then $-x\eq b\pmod 4$.)
Thus the numbers in (\ref{2.6}) are all integers.
\medskip

{\it Case} 2. $2\|a$ and $2\mid b$.

Write $a=2a_0$ and $b=2b_0$ with $a_0,b_0\in\Z$ and $2\nmid a_0$.  Since $2a_0-b_0^2\equiv 1,2\pmod{4}$, we have $2a_0-b_0^2=x_0^2+y_0^2+z_0^2$ for some $x_0,y_0,z_0\in\Z$.
Without loss of generality, we may assume that $x_0\equiv b_0\pmod{2}$ and $y_0\equiv z_0\pmod{2}$.
Set $x=2x_0$, $y=2y_0$ and $z=2z_0$. Then $4a-b^2=4(2a_0-b_0^2)=x^2+y^2+z^2$, and the numbers in (\ref{2.6}) are all integers.

 In either case, there are $x,y,z\in\Z$ for which $4a-b^2=x^2+y^2+z^2$ and $s,t,u,v\in\Z$, where $s,t,u,v$ are as in (\ref{2.6}).
 Obviously, $s+t+u+v=b$ and
\begin{align*}&s^2+t^2+u^2+v^2
\\=&2\l(\f{s+t}2\r)^2+2\l(\f{s-t}2\r)^2+2\l(\f{u+v}2\r)^2+2\l(\f{u-v}2\r)^2
\\=&2\l(\f{b+x}4\r)^2+2\l(\f{y+z}4\r)^2+2\l(\f{b-x}4\r)^2+2\l(\f{y-z}4\r)^2
\\=&\f{b^2+x^2+y^2+z^2}4=a.
\end{align*}
In view of Lemma \ref{Lem2.3} and the second inequality in (\ref{ineq}), we have
\[(|x|+|y|+|z|)^2\ls3(x^2+y^2+z^2)=3(4a-b^2)<(b+4)^2.\]
Therefore
\[s,t,u,v\gs\frac{b-|x|-|y|-|z|}4>-1\]
and hence $s,t,u,v\in\N$.  \qed

Now we need one more lemma which is well known (cf. \cite[p.\,59]{B}).

\begin{lemma}\label{Lem2.5} For any $n\in\Z^+$, we have
\begin{equation}\label{2.7}r_4(n)=8\sum_{d\mid n\atop 4\nmid d} d\end{equation}
where
$$r_4(n):=|\{(w,x,y,z)\in\Z^4:\ w^2+x^2+y^2+z^2=n\}|.$$
\end{lemma}

\medskip
\noindent{\it Proof of Theorem 1.1}. (i) Let $I_1=[\al,\beta]$ be the interval given by (\ref{2.1}).
As $N\gs 7\times 2^2m^3$, by Lemma 2.1 the length of the interval $I_1$ is greater than $2m$.
Choose $b_0\in\{\lceil\alpha\rceil+r:\ r=0,\ldots,m-1\}$ with $b_0\eq N\pmod m$. Then $b_1=b_0+m\ls\lceil\alpha\rceil+2m-1<\al+2m<\beta.$
Thus both $b_0$ and $b_1$ lie in $I_1$.
Note that
$$\f2m(N-b_1)+b_1-\l(\f2m(N-b_0)+b_0\r)= m-2\eq 2\pmod 4.$$
So, for some $b\in\{b_0,b_1\}$ and $a=\f2m(N-b)+b$ we have $2\nmid ab$, or $2\| a$ and $2\mid b$.
Obviously,
$$b\gs\min I_1>0,\ a=\f 2mN+\l(1-\f2m\r)b>0,\ \mbox{and}\ N=\f m2(a-b)+b.$$
Applying Lemmas 2.2 and 2.4, we see that there are $s,t,u,v\in\N$ satisfying (\ref{2.4}).
Therefore,
\begin{align*}N=&\f m2(a-b)+b
\\=&\f m2(s^2-s+t^2-t+u^2-u+v^2-v)+s+t+u+v
\\=&p_{m+2}(s)+p_{m+2}(t)+p_{m+2}(u)+p_{m+2}(v).
\end{align*}

(ii) Write $m=4l$ with $l\in\Z^+$. Let $\varphi$ denote Euler's totient function. We want to show that none of the positive integers
$$4l^2\times\f{4^{k\varphi(2l+1)}-1}{2l+1}\ \ (k=1,2,\ldots)$$
can be written as $\sum_{i=1}^4 p_{m+4}(x_i)$ with $x_1,x_2,x_3,x_4\in\N$.

Suppose that for some $n\in\Z^+$ divisible by $\varphi(2l+1)$ there are $w,x,y,z\in\N$ such that
\begin{align*}4l^2\times\f{4^n-1}{2l+1}=&p_{m+4}(w)+p_{m+4}(x)+p_{m+4}(y)+p_{m+4}(z)
\\=&\f {4l+2}2(w^2+x^2+y^2+z^2-w-x-y-z)+w+x+y+z.
\end{align*}
Then
$$4^{n+1}l^2=((2l+1)w-l)^2+((2l+1)x-l)^2+((2l+1)y-l)^2+((2l+1)z-l)^2.$$
As $r_4(4^{n+1}l^2)=r_4(4l^2)$ by Lemma 2.5, there are $w_0,x_0,y_0,z_0\in \Z$ with $w_0^2+x_0^2+y_0^2+z_0^2=4l^2$ such that
\begin{align*}&(2l+1)w-l=2^nw_0,\ (2l+1)x-l=2^nx_0,
\\ &(2l+1)y-l=2^ny_0, \ (2l+1)z-l=2^nz_0.
\end{align*}
Since $2^n\equiv1\pmod{2l+1}$ by Euler's theorem, we have
$$w_0\eq x_0\eq y_0\eq z_0\equiv -l\pmod{2l+1}.$$
As $w_0^2+x_0^2+y_0^2+z_0^2=4l^2$, we must have $w_0=x_0=y_0=z_0=-l$.
Thus $(2l+1)w=2^nw_0+l=l(1-2^n)<0$, which contradicts $w\in\N$. \qed

\section{Proof of Theorem 1.2}
\setcounter{lemma}{0}
\setcounter{theorem}{0}
\setcounter{corollary}{0}
\setcounter{remark}{0}
\setcounter{equation}{0}
\setcounter{conjecture}{0}

\begin{lemma}\label{Lem3.1} Let $l,m,N\in\Z^+$ with $N\gs11lm^2(lm+1)$. Then the length of the interval
\begin{equation}\label{3.1} I_2=\left[ \frac32+\sqrt{\frac{10N}m-3},\ 1+\sqrt\frac {12N}m \right]\end{equation}
is greater than $lm$.
\end{lemma}
\Proof. The length $L_2$ of the interval $I_2$ is $\sqrt{12x}-\sqrt{10x-3}-1/2$, where $x=N/m\gs 11lm(lm+1)$. Let
$l_0=lm+1/2$. Then
\begin{align*}L_2>lm&\iff \sqrt{12x}>\sqrt{10x-3}+l_0
\\&\iff 2x+3-l_0^2+3>2l_0\sqrt{10x-3}
\\&\iff 4x(x-11l_0^2+3)+(l_0^2-3)^2+12l_0^2>0.
\end{align*}
As
$$x\gs 11lm(lm+1)>11l^2m^2+11lm+\f{11}4-3=11l_0^2-3,$$
we have $L_2>lm$ as desired. \qed

\begin{lemma}\label{Lem3.2} Let $a,b,m,N\in\Z^+$ with $m\gs 3$ and
$$N=\frac m2(a-b)+b\gs \f35m.$$
Suppose that $b$ belongs to the interval $I_2$ given by $(\ref{3.1})$.
Then
\begin{equation}\label{3.2} b^2<6a\ \ \mbox{and}\ \ 5a<b^2+2b+6. \end{equation}
\end{lemma}
\Proof. Note that
\[a=\left(1-\frac2m\right)b+\frac {2N}m.\]
Thus
\begin{align*}b^2-6a=&b^2-6\left(1-\frac2m\right)b-\frac{12N}m
\\=&\l(b-3\l(1-\f2m\r)\r)^2-9\l(1-\f2m\r)^2-\f{12N}m.
\end{align*}
As $b\in I_2$, we have
$$ b\ls 1+\sqrt\frac {12N}m<3\left(1-\frac2m\right)+\sqrt{9\left(1-\frac2m\right)^2+\frac {12N}m}$$
and hence $b^2-6a<0$.
On the other hand,
\begin{align*}b^2+2b+6-5a=&b^2-\left(3-\frac{10}m\right)b+\left(6-\frac{10N}m\right)
\\=&\l(b-\l(\f 32-\f 5m\r)\r)^2-\l(\f 32-\f 5m\r)^2+\left(6-\frac{10N}m\right)
\\>&0
\end{align*}
since $(3/2-5/m)^2\ls 9/4<3$ and
$$b\gs\frac32+\sqrt{\frac{10N}m-3}>\frac32-\frac5m+\sqrt{\left(\frac32-\frac5m\right)^2-6+\frac{10N}m}.$$
This proves (\ref{3.2}). \qed

\begin{lemma}\label{Lem3.3} Let $a$ and $b$ be positive integers with $a\eq b\pmod 2$ satisfying $(\ref{3.2})$.
Suppose that $2\nmid a$ or $4\|a$, and that $3\mid a$ or $3\nmid b$.
Then there exist $s,t,u,v\in\N$ such that
\begin{equation}\label{3.3}a=s^2+t^2+2u^2+2v^2
\ \mbox{and}\ b=s+t+2u+2v.\end{equation}
\end{lemma}
\Proof. If $n\in\N$ is not of the form $4^k(8l+7)$ with $k,l\in\N$, then by (\ref{2.5}) there are $x,u,v\in\Z$ with $u\eq v\pmod 2$ such that
$$n=x^2+u^2+v^2=x^2+2\l(\f{u-v}2\r)^2+2\l(\f{u+v}2\r)^2.$$
We claim that there are $x,y,z\in\Z$ with $6a-b^2=x^2+2y^2+2z^2$ such that all the numbers
\begin{equation}\label{3.4}\begin{aligned} s=\frac{b+x+2y+2z}6,&\ \ t=\frac{b-x-2y+2z}6,
\\u=\frac{b-x+y-z}6,&\ \ v=\frac{b+x-y-z}6
\end{aligned}
\end{equation}
are integral.
\medskip

{\it Case} 1. $3\nmid b$.

If $a\eq b\eq1\pmod2$, then $6a-b^2\equiv1\pmod{4}$. When $4\|a$ and $2\mid b$, we have $6a-b^2\eq4,8\pmod{16}$.
 Thus $6a-b^2=x^2+2y^2+2z^2$ for some $x,y,z\in\Z$. Clearly, $x\eq b\pmod2$ and $y\equiv z\pmod{2}$.
 Note that
 $$x^2+2y^2+2z^2\equiv 6a-b^2\equiv2\pmod3.$$
 As $x^2\not\eq2\pmod 3$, we have $3\nmid y$ or $3\nmid z$.
 Without loss of generality, we assume that $3\nmid z$ and moreover $z\eq b\pmod 3$.
 (If $z\eq-b\pmod 3$ then $-z\eq b\pmod3$.) As $x^2-y^2\eq x^2+2y^2\eq0\pmod 3$,
 without loss of generality we may assume that $x\eq y\pmod 3$.
 Now it is easy to see that all the numbers in (\ref{3.4}) are integers.
\medskip

{\it Case} 2. $3\mid a$ and $3\mid b$.

In this case,  $a=3a_0$ and $b=3b_0$ for some $a_0,b_0\in\Z$.
If $a\eq b\eq1\pmod2$, then $2a_0-b_0^2\equiv1\pmod{4}$.
When $4\|a$ and $2\mid b$, we have $2a_0-b_0^2\eq4,8\pmod{16}$. Thus
$2a_0-b_0^2=x^2_0+2y_0^2+2z_0^2$ for some $x_0,y_0,z_0\in\Z$. It follows that $x_0\eq b_0\pmod2$ and $y_0\equiv z_0\pmod{2}$ since $a_0\eq b_0\pmod2$.
Set $x=3x_0$, $y=3y_0$ and $z=3z_0$. Then $6a-b^2=9(2a_0-b_0^2)=x^2+2y^2+2z^2$ and all the numbers in (\ref{3.4}) are integers.

By the above, in either case, there are $x,y,z\in\Z$ for which $x^2+2y^2+2z^2=6a-b^2$ and $s,t,u,v\in\Z$, where $s,t,u,v$ are as in (\ref{3.4}).
Observe that
$$s+t+2(u+v)=\f{b+2z}3+2\times\f{b-z}3=b$$ and also
\begin{align*}&s^2+t^2+2u^2+2v^2
\\=&2\l(\f{s+t}2\r)^2+2\l(\f{s-t}2\r)^2+(u+v)^2+(u-v)^2
\\=&2\l(\f{b+2z}6\r)^2+2\l(\f{x+2y}6\r)^2+\l(\f{b-z}3\r)^2+\l(\f{x-y}3\r)^2
\\=&\f{b^2+x^2+2y^2+2z^2}6=a.
\end{align*}

In view of Lemma \ref{Lem2.3} and (\ref{3.2}),
\[(|x|+2|y|+2|z|)^2\ls5(x^2+2y^2+2z^2)=5(6a-b^2)<(b+6)^2\]
and hence
\[b-|x|-2|y|-2|z|>-6.\]
So we have
\[s,t,u,v\gs\frac{b-|x|-2|y|-2|z|}6>-1,\]
and hence $s,t,u,v\in\N$.

In view of the above, we have completed the proof of Lemma 3.3. \qed

\medskip
\noindent{\it Remark} 3.1. For $s,t,u,v$ given in (\ref{3.4}), the identity
\begin{align*}6(s^2+t^2+2u^2+2v^2)=&b^2+x^2+2y^2+2z^2
\\=&(s+t+2u+2v)^2+(s-t-2u+2v)^2
\\&+2(s-t+u-v)^2+2(s+t-u-v)^2
\end{align*}
is a special case of our following general identity
\begin{equation}\label{identity}\begin{aligned}&(a+b)(c+d)(w^2+abx^2+cdy^2+abcdz^2)
\\=&ac(w+bx+dy+bdz)^2+ad(w+bx-cy-bcz)^2
\\&+bc(w-ax+dy-adz)^2+bd(w-ax-cy+acz)^2.
\end{aligned}\end{equation}
We have also found another similar identity:
\begin{equation}\label{3.6}\begin{aligned} &(3b+4)(w^2+2x^2+(b+1)y^2+2bz^2)
\\=&(w+2x+(b+1)y+2bz)^2+2(w-(b+1)y+bz)^2
\\&+(b+1)(w-2x+y)^2+2b(w+x-2z)^2.
\end{aligned}\end{equation}

\medskip
\noindent{\it Proof of Theorem 1.2}. (i) Let $I_2=[\alpha,\beta]$ be the interval given by (\ref{3.1}). As $N\gs 1628m^3=(12+1/3)m\times132m^2\gs11m^2\times 12(12m+1)$,
by Lemma 3.1 the length of the interval $I_2$ is greater than $12m$. We distinguish two cases to construct integers $b\in I_2$ and $a\eq b\pmod 2$
for which $N=\f m2(a-b)+b$, and $2\nmid a$ or $4\|a$, and $3\mid a$ or $3\nmid b$.
\medskip

{\it Case} 1. $3\nmid m$ or $3\nmid N$.

Choose $b_0\eq N\pmod m$ with $b_0\in\{\lceil\al\rceil+r:\ r=0,\ldots,m-1\}$, and let $b_j=b_0+jm$ for $j=1,\ldots,7$.
Since $\lceil \al\rceil+8m-1<\al+8m<\beta$, we have $b_i\in I_2$ for all $i=0,\ldots,7$.

If $2\nmid m$, then we choose $i\in\{0,1\}$ with $b_i$ odd. When $4\mid m$, we may choose $i\in\{0,1,2,3\}$ with
$$a_i:=\f 2m(N-b_i)+b_i\eq4\pmod 8$$
since $a_i-a_0=-2i+im=2i(m/2-1)$ with $m/2-1$ odd.
If $3\mid m$ and $3\nmid N$, then $b=b_i\eq N\not\eq0\pmod3$.
When $3\nmid m$, we choose $j\in\{i,i+4\}$ with $b=b_j\not\eq0\pmod 3$, and note that
$$a_{i+4}-a_i=4m-8\eq\begin{cases}0\pmod2&\mbox{if}\ 2\nmid m,\\0\pmod8&\mbox{if}\ 4\mid m.\end{cases}$$

As $N\eq b\pmod m$, we see that $a=2(N-b)/m+b$ is an integer with $a\eq b\pmod2$.
By our choice, $2\nmid b$ if $2\nmid m$, and $a\eq 4\pmod 8$ if $4\mid m$.  Note also that $3\nmid b$.
\medskip

{\it Case} 2. $m\eq N\eq0\pmod 3$.

Choose $b_0\in\{\lceil\al\rceil+r:\ r=0,1,\ldots,3m-1\}$ with $b_0\eq N\pmod{3m}$. If $2\nmid m$, then we choose $b\in\{b_0,b_0+3m\}$ with $b$ odd and hence
$a=\f2m(N-b)+b\eq b\eq1\pmod2$.
When $4\mid m$, we may choose $b\in\{b_0+3jm:\ j=0,1,2,3\}$ with $a=\f2m(N-b)+b\eq4\pmod 8$, for
$$\f2m(N-b_0-3jm)+b_0+3jm-\l(\f2m(N-b_0)+b_0\r)=6j\l(\f m2-1\r)$$
with $m/2-1$ odd.
Note that
$$\al\ls b_0<b_0+9m\ls\lceil\al\rceil+3m-1+9m<\al+12m<\beta$$
and hence $b\in I_2$. Obviously, $a\eq b\eq N\eq0\pmod 3.$

Now we have constructed positive integers $b\in I_2$ and $a\eq b\pmod2$ with $N=\f m2(a-b)+b$
such that $2\nmid a$ or $4\|a$, and $3\mid a$ or $3\nmid b$. So (3.2) holds by Lemma 3.2. In view of Lemma 3.3, (3.3) holds for some $s,t,u,v\in\N$.
Therefore,
\begin{align*}N=&\f m2(s^2+t^2+2u^2+2v^2-s-t-2u-2v)+s+t+2u+2v
\\=&m\bi s2+m\bi t2+2m\bi u2+2m\bi v2+s+t+2u+2v
\\=&p_{m+2}(s)+p_{m+2}(t)+2p_{m+2}(u)+2p_{m+2}(v).
\end{align*}
This proves part (i) of Theorem 1.2.

(ii) Now assume that $m=2l$ with $l\in\Z^+$ odd. We want to show that none of the  positive integers
$$(l-1)^2\times \f{4^{k\varphi(l)}-1}l\ \ (k=1,2,3,\ldots)$$
can be written as $p_{m+2}(w)+p_{m+2}(x)+2p_{m+2}(y)+2p_{m+2}(z)$ with $w,x,y,z\in\N$.

Let $n\in\Z^+$ be a multiple of $\varphi(l)$. Then $2^n\eq1\pmod l$ by Euler's theorem. Suppose that there are $w,x,y,z\in\N$ for which
\begin{align*}&(l-1)^2\times\f{4^n-1}l
\\=&p_{m+2}(w)+p_{m+2}(x)+2p_{m+2}(y)+2p_{m+2}(z)
\\=&\f{2l}2(w^2+x^2+2y^2+2z^2-w-x-2y-2z)+w+x+2y+2z.
\end{align*}
Then we have
\begin{align*}&4^{n+1}(l-1)^2
\\=&(2lw-(l-1))^2+(2lx-(l-1))^2+2(2ly-(l-1))^2+2(2lz-(l-1))^2
\\=&(2lw-(l-1))^2+(2lx-(l-1))^2+(2l(y+z)-2(l-1))^2+(2l(y-z))^2
\end{align*}
and hence
$$4^n(l-1)^2=\l(lw-\f{l-1}2\r)^2+\l(lx-\f{l-1}2\r)^2+(l(y+z-1)+1)^2+(l(y-z))^2.$$
As $4\mid(l-1)^2$, by Lemma \ref{Lem2.5} we have $r_4(4^n(l-1)^2)=r_4((l-1)^2)$.
So there are $w_0,x_0,y_0,z_0\in\Z$ with
\begin{equation}\label{3.7}w_0^2+x_0^2+y_0^2+z_0^2=(l-1)^2\end{equation}
such that
$$lw-\f{l-1}2=2^nw_0,\ lx-\f{l-1}2=2^nx_0,\ l(y+z-1)+1=2^ny_0,\ l(y-z)=2^nz_0.$$
As $2^n\eq1\pmod l$, we see that
$$w_0\eq\f{1-l}2\pmod l,\ x_0\eq\f{1-l}2\pmod l,\ y_0\eq1\pmod l,\ z_0\eq0\pmod l.$$
Observe that $l=m/2\gs2$ and hence $w_0x_0\not=0$. Thus $y_0^2+z_0^2\ls(l-1)^2-2$.
Since $y_0\eq1\pmod l$ and $z_0\eq0\pmod l$, we must have $y_0=1$ and $z_0=0$.
Now (\ref{3.7}) yields $w_0^2+x_0^2=(l-1)^2-1=l^2-2l$. As
$w_0\eq x_0\eq (1-l)/2\pmod l$, we must have $\{w_0,x_0\}\se\{(1-l)/2,(1+l)/2\}$.
It is easy to verify directly that none of the numbers
$$\l(\f{1-l}2\r)^2+\l(\f{1+l}2\r)^2,\ 2\l(\f{l-1}2\r)^2,\ 2\l(\f{1+l}2\r)^2$$
is equal to $(l-1)^2-1=l^2-2l$. This contradiction concludes the proof of
Theorem 1.2(ii). \qed

\section{Proof of Theorem 1.3}
\setcounter{lemma}{0}
\setcounter{theorem}{0}
\setcounter{corollary}{0}
\setcounter{remark}{0}
\setcounter{equation}{0}
\setcounter{conjecture}{0}

\begin{lemma}\label{Lem4.1} Let $a$ and $b$ be positive integers satisfying $(3.2)$ for which $a\eq b\pmod 2$, and $a\eq3\pmod9$ or $3\nmid b$.
Then there exist $s,t,u,v\in\N$ such that
\begin{equation}\label{4.1}a=s^2+t^2+u^2+3v^2\ \ \mbox{and}\ \ b=s+t+u+3v.\end{equation}
\end{lemma}
\Proof. It is known that (cf. Dickson \cite[pp.112-113]{D39})
\begin{equation}\label{4.2}\{x^2+y^2+3z^2:\ x,y,z\in\Z\}=\N\sm\{9^k(9l+6):\ k,l\in\N\}.
\end{equation}

If $3\nmid b$, then $6a-b^2\equiv2\pmod{3}$. If $a\eq3\pmod9$ and $3\mid b$, then $6a-b^2\equiv\pm9\pmod{27}$.
By (\ref{4.2}), $6a-b^2=x^2+y^2+3z^2$ for some $x,y,z\in\Z$.
Clearly, $x^2+y^2\equiv2b^2\pmod{3}$. Without loss of generality, we may assume that $x\equiv y\equiv b\pmod3$. (If $x\eq-b\pmod 3$ then $-x\eq b\pmod3$.)
If $a$ and $b$ are both odd, then $x^2+y^2+3z^2=6a-b^2\equiv1\pmod{4}$, and hence one of $x$ and $y$ is odd.
If $a$ and $b$ are both even, then $x^2+y^2+3z^2=6a-b^2\equiv0\pmod{4}$, and hence one of $x$ and $y$ is even.
Without loss of generality, we may assume that $x\eq a\eq b\pmod 2$ and $y\equiv z\pmod2$.
Thus all the numbers
\begin{equation}\label{4.3}s=\frac{b+x+y+3z}6,\ t=\frac{b+ x+y-3z}6, \ u=\frac{b+x-2y}6,
\ v=\frac{b-x}6
\end{equation}
are integers.
Observe that
$$s+t+u+3v=\f{b+x}2+3\times\f{b-x}6=b$$
and
\begin{align*}&s^2+t^2+u^2+3v^2
\\=&2\l(\f{s+t}2\r)^2+2\l(\f{s-t}2\r)^2+u^2+3v^2
\\=&2\l(\f{b+x+y}6\r)^2+2\l(\f z2\r)^2+\l(\f{b+x-2y}6\r)^2+3\l(\f{b-x}6\r)^2
\\=&\f{b^2+x^2+y^2+3z^2}6=a.
\end{align*}

In view of Lemma \ref{Lem2.3} and (3.2),
\[(|x|+|y|+3|z|)^2\ls5(x^2+y^2+3z^2)=5(6a-b^2)<(b+6)^2\]
and hence
\[b-|x|-|y|-3|z|>-6.\]
So we have
\[s,t,u,v\gs\frac{b-|x|-|y|-3|z|}6>-1\]
and hence $s,t,u,v\in\N$. This concludes the proof. \qed

\medskip
\noindent{\it Proof of Theorem 1.3}.
As $$N\gs924m^3=99m^3\l(9+\f13\r)\gs99m^3\l(9+\f1m\r)= 99m^2(9m+1),$$
 the length of the interval $I_2=[\al,\beta]$ defined in Lemma 3.1 is greater than $9m$.

Let $b_0\in\{\lceil\al\rceil+r:\ r=0,1,\ldots,m-1\}$ with $b_0\eq N\pmod m$.
If $3\nmid m$ or $3\nmid N$, then we may choose $b\in\{b_0,b_0+m\}$ with $b\not\eq0\pmod3$.
When $3\mid m$ and $3\mid N$, we let $c_0\in\{\lceil\al\rceil+r:\ r=0,\ldots,3m-1\}$ with $c_0\eq N\pmod{3m}$,
and set $b=c_0+j3m$ with $j\in\{0,1,2\}$ such that
$$\f 2m(N-c_0)+c_0-6j\eq3\pmod 9.$$ Note that $b\in I_2$ since
$$\al\ls b\ls \lceil\al\rceil+3m-1+6m<\al+9m<\beta.$$

Let
 $$a=\frac2m(N-b)+b,\ \ \mbox{i.e.,}\ N=\f m2(a-b)+b.$$
Then
$$a=\f2mN+\l(1-\f2m\r)b>0\quad\mbox{and}\quad a\eq b\pmod 2.$$
If $3\mid b$, then $3\mid m$ and
$$a=\f2m(N-b)+b=\f2m(N-c_0-3jm)+c_0+3jm\eq\f2m(N-c_0)-6j\eq 3\pmod 9.$$

By Lemmas 3.2 and 4.1, there are $s,t,u,v\in\N$ satisfying (4.1). Therefore,
\begin{align*}N=&\f m2(s^2+t^2+u^2+3v^2-s-t-u-3v)+s+t+u+3v
\\=&p_{m+2}(s)+p_{m+2}(t)+p_{m+2}(u)+3p_{m+2}(v).
\end{align*}
This completes the proof of Theorem 1.3. \qed

\section{Proof of Theorem 1.4}
\setcounter{lemma}{0}
\setcounter{theorem}{0}
\setcounter{corollary}{0}
\setcounter{remark}{0}
\setcounter{equation}{0}
\setcounter{conjecture}{0}

\begin{lemma}\label{Lem5.1} Let $l,m,N\in\Z^+$ with $lm\gs 20$ and $N\gs3lm^2(5lm+12)$. Then the length of the interval
\begin{equation}\label{5.1}I_3= \left[ \frac52+\sqrt{\frac{14N}m-1},\ \frac43+4\sqrt\frac Nm \right]\end{equation}
is greater than $lm$.
\end{lemma}
\Proof. The length of the interval $I_3$ is $4\sqrt x-\sqrt{14x-1}-7/6$, where $x=N/m$.
Set $l_0=lm+7/6$. Then
\begin{align*} &4\sqrt x-\sqrt{14x-1}-\f76>lm
\\\iff&4\sqrt x>\sqrt{14x-1}+l_0
\\\iff&2x+1-l_0^2>2l_0\sqrt{14x-1}
\\\iff&4x(x+1-15l_0^2)+(l_0^2-1)^2+4l_0^2>0.
\end{align*}
Note that
$$x\gs 15l^2m^2+36lm>15l^2m^2+35lm+\f{245}{12}-1=15l_0^2-1.$$
So the desired result follows. \qed

\begin{lemma}\label{Lem5.2} Let $a,b,m,N$ be positive integers with $m\gs3$,
\begin{equation*}N=\frac m2(a-b)+b\gs\f47m\end{equation*}
and $b\in I_3$, where $I_3$ is the interval given by $(\ref{5.1})$.
Then
\begin{equation}\label{5.2} b^2<8a \ \ \mbox{and}\ \ 7a<b^2+2b+8. \end{equation}
\end{lemma}
\Proof. Clearly,
\[a=\left(1-\frac2m\right)b+\frac {2N}m.\]
Thus
\begin{align*}b^2-8a=&b^2-8\left(1-\frac2m\right)b-\frac{16N}m
\\=&\l(b-4\l(1-\f2m\r)\r)^2-16\l(\l(1-\f2m\r)^2+\f Nm\r).
\end{align*}
As $b\in I_3$, we have
$$b\ls\frac43+4\sqrt\frac Nm<4\left(1-\frac2m\right)+4\sqrt{\left(1-\frac2m\right)^2+\frac Nm}$$
and hence $b^2-8a<0$.
On the other hand,
\begin{align*}&b^2+2b+8-7a
\\=&b^2-\left(5-\frac{14}m\right)b+\left(8-\frac{14N}m\right)
\\=&\l(b-\left(\f 52-\f 7m\r)\r)^2-\l(\l(\f 52-\f 7m\r)^2+\f{14N}m-8\r)>0
\end{align*}
since $(5/2-7/m)^2\ls 25/4<7$ and
$$b\gs\frac52+\sqrt{\frac{14N}m-1}>\frac52-\frac7m+\sqrt{\left(\frac52-\frac7m\right)^2+\frac{14N}m-8}.$$
Therefore (\ref{5.2}) holds. \qed

\begin{lemma}\label{Lem5.3} Let $a$ and $b$ be positive integers satisfying $(\ref{5.2})$.
Then there are $s,t,u,v\in\N$ such that
\begin{equation}\label{5.3}a=s^2+t^2+2u^2+4v^2\ \ \mbox{and}\ \ b=s+t+2u+4v,\end{equation}
under one of the following conditions {\rm (i)-(iii):}

{\rm (i)} $2\nmid ab$.

{\rm (ii)} $2\mid a$ and $2\|b$.

{\rm (iii)} $4\mid a$ and $4\|b$, or $a\equiv b+4\pmod{16}$ and $8\mid b$.
\end{lemma}
\Proof. It is known (cf. \cite[pp.\,112-113]{D39}) that
\begin{equation}\label{5.4}\{x^2+2y^2+4z^2:\ x,y,z\in\Z\}=\N\sm\{4^k(16l+14):\ k,l\in\N\}.\end{equation}
We claim that if one of (i)-(iii) holds then $8a-b^2=x^2+2y^2+4z^2$ for some $x,y,z\in\Z$ such that all the numbers
\begin{equation}\label{5.5}u=\frac{b+x-2y}8,\ v=\frac{b-x}8,
\ s=u+\frac{y+z}2,\ t=u+\frac{y-z}2
\end{equation}
are integers.
\medskip

{\it Case} 1. $2\nmid ab$.

Since $8a-b^2\equiv-1\pmod{8}$, we have $8a-b^2=x^2+2y^2+4z^2$ for some $x,y,z\in\Z$. As $x^2+2y^2\equiv-1\pmod{4}$, we have $2\nmid xy$.
Since $4z^2\eq -b^2-x^2-2y^2\eq-1-1-2\pmod 8$, we also have $2\nmid z$. Note that
$$x^2=8a-b^2-2y^2-4z^2\eq 8-b^2-2-4=2-b^2\eq b^2\pmod{16}$$
and hence $x\eq\pm b\pmod8$. Without loss of generality, we may assume that $x\eq b\pmod 8$ and $y\eq b\pmod 4$. (If $y\eq-b\pmod4$ then $-y\eq b\pmod4$.)
Thus all the four numbers in (\ref{5.5}) are integers.
\medskip

{\it Case} 2. $2\mid a$ and $2\|b$.

Write $a=2a_0$ and $b=2b_0$ with $a_0,b_0\in\Z$ and $2\nmid b_0$.
Since $4a_0-b_0^2\equiv3\pmod{4}$, by (\ref{5.4}) we have $4a_0-b_0^2=x_0^2+2y_0^2+4z_0^2$ for some $x_0,y_0,z_0\in\Z$.
As $x_0^2+2y_0^2\equiv3\pmod{4}$, both $x_0$ and $y_0$ are odd. Without loss of generality, we may assume that $x_0\eq y_0\eq b_0\pmod 4$.
Set $x=2x_0$, $y=2y_0$ and $z=2z_0$. Then
$$8a-b^2=4(4a_0-b_0^2)=x^2+2y^2+4z^2$$
and all the numbers in (\ref{5.5}) are integers.
\medskip

{\it Case} 3. $4\mid a$ and $4\|b$, or $a\eq b+4\pmod{16}$ and $8\mid b$.

Write $a=4a_0$ and $b=4b_0$ with $a_0,b_0\in\Z$. Then
$$2a_0-b_0^2\eq\begin{cases}1\pmod{2}&\mbox{if}\ 4\|b\ (\mbox{i.e.,}\ 2\nmid b_0),
\\2(b_0+1)-b_0^2\eq2\pmod 8&\mbox{if}\ a\eq b+4\pmod{16}\ \mbox{and}\ 8\mid b.
\end{cases}$$
Thus, by (\ref{5.4}) we have $2a_0-b_0^2=x_0^2+2y_0^2+4z_0^2$ for some $x_0,y_0,z_0\in\Z$. Obviously $x_0\eq b_0\pmod 2$.
Set $x=4x_0$, $y=4y_0$ and $z=4z_0$. Then all the numbers in (\ref{5.5}) are integers.

Now assume that one of the conditions (i)-(iii) holds. By the claim we have proved, there are $x,y,z\in\Z$ such that
$8a-b^2=x^2+2y^2+4z^2$ and also $s,t,u,v\in\Z$, where $s,t,u,v$ are given by (\ref{5.5}).
Clearly,
$$s+t+2u+4v=y+2u+2u+4v=\f{b+x}2+\f{b-x}2=b$$
and
\begin{align*}&s^2+t^2+2u^2+4v^2
\\=&2\l(u+\f y2\r)^2+2\l(\f z2\r)^2+2\l(\f{b+x-2y}8\r)^2+4\l(\f{b-x}8\r)^2
\\=&2\l(\f{b+x+2y}8\r)^2+2\l(\f{b+x-2y}8\r)^2+\f{z^2}2+\l(\f{b-x}4\r)^2
\\=&\l(\f{b+x}4\r)^2+\l(\f y2\r)^2+\f{z^2}2+\l(\f{b-x}4\r)^2
\\=&\f{b^2+x^2+2y^2+4z^2}8=a.
\end{align*}

In view of Lemma \ref{Lem2.3} and (\ref{5.2}),
\[(|x|+2|y|+4|z|)^2\ls7(x^2+2y^2+4z^2)=7(8a-b^2)<(b+8)^2\]
and hence
\[b-|x|-2|y|-4|z|>-8.\]
So we have
\[u,v,s,t\gs\frac{b-|x|-2|y|-4|z|}8>-1\]
and hence $s,t,u,v\in\N$.

In view of the above, we have completed the proof of Lemma 5.3. \qed

\medskip

\noindent{\it Proof of Theorem 1.4}. As
$$N\gs 96m^2\times 11m \gs 96m^2(10m+3)=24m^2(40m+12),$$
applying Lemma 5.1 with $l=8$ we find that
the interval $I_3=[\al,\beta]$ given by (\ref{5.1}) has length greater than $8m$.
\medskip

{\it Case} 1. $4\nmid m$ or $8\nmid N$.

Choose $b_0\in\{\lceil\al\rceil+r:\ r=0,\ldots,m-1\}$ with $N\equiv b_0\pmod{m}$.
Set $b_1=b_0+m$ . Then
$$\al\ls b_0<b_1\ls \lceil\al\rceil+2m-1<\al+8m<\beta$$
and thus $b_j\in I_3$ for each $j=0,1$. Note that
$$a_j:=\f2m(N-b_j)+b_j=\f2m(N-b_0)+b_0+(m-2)j.$$

If $2\nmid m$, then $a_j\eq b_j\eq1\pmod2$ for some $j\in\{0,1\}$.
If $2\mid m$ and $2\nmid N$, then $a_0\eq b_0\eq1\pmod 2$.

If $2\|m$ and $2\mid N$, then for some $j\in\{0,1\}$ we have $b_j\eq2\pmod 4$ and $2\mid a_j$.
When $4|m$ and $2\|N$, we have $b_0\eq2\pmod 4$ and $a_0\eq b_0\eq0\pmod 2$.
If $4\|m$ and $4\| N$, then $4\mid b_0$, and for some $j\in\{0,1\}$ we have
$b_j\eq4\pmod 8$ and $a_j\eq b_j\eq0\pmod 4$.
When $8\mid m$ and $4\|N$, we have $b_0\eq4\pmod 8$ and $a_0\eq b_0\pmod2$,
hence for some $j\in\{0,1\}$ we have $a_j\eq0\pmod 4$ and $b_j\eq b_0\eq4\pmod 8$.
\medskip

{\it Case} 2. $4\mid m$ and $8\mid N$.

Choose $b\in\{\lceil\al\rceil+r:\ r=0,\ldots,8m-1\}$ such that $b\eq N-2m\pmod{8m}$.
Since $\al\ls b\ls\lceil\al\rceil+8m-1<\al+8m<\beta$, we have $b\in I_3$.
Clearly, $8\mid b$ since $8\mid N$ and $4\mid m$. Note that
$$\f2m(N-b)+b\eq 4+b\pmod 8.$$

By the above, in either case we can always find $b\in I_3$ and $a\in\Z$ for which one of (i)-(iii) in Lemma 5.3 holds and also
$$a=\f2m(N-b)+b,\ \ \mbox{i.e.,}\ N=\f m2(a-b)+b.$$
By Lemmas 5.2 and 5.3, there are $s,t,u,v\in\N$ satisfying (\ref{5.3}). Therefore,
\begin{align*}N=&\f m2(s^2+t^2+2u^2+4v^2-s-t-2u-4v)+s+t+2u+4v
\\=&p_{m+2}(s)+p_{m+2}(t)+2p_{m+2}(u)+4p_{m+2}(v).
\end{align*}

In view of the above, we have finished the proof of Theorem 1.4. \qed

\subsection*{Acknowledgements}
This research was supported by the National Natural Science Foundation of China (grant no. 11571162).

\end{document}